\documentclass[12pt]{amsart}
\usepackage{amssymb}

\newtheorem{Thm}{Theorem}[section]

\newtheorem{Lem}[Thm]{Lemma}
\newtheorem{Prop}[Thm]{Proposition}

\theoremstyle{definition}
\newtheorem{Def}[Thm]{Definition}

\theoremstyle{remark}

\def\ldots{\mathinner{\ldotp\ldotp\ldotp}}
\def\ldots{\mathinner{\cdotp\cdotp\cdotp}}
\def \tree{\omega^{<\omega}}

\def \eps{\varepsilon}
\def\Cdb{\mathbb C}
\def\Ndb{\mathbb N}

\def\tnorm#1{\vert\vert\vert#1\vert\vert\vert} 

\begin{document}

\title{$L^p-$maximal regularity on Banach spaces with a Schauder basis}
\author{N. J. Kalton}
\address{Department of Mathematics \\
University of Missouri-Columbia \\ Columbia, MO 65211 }
\thanks{The first author was partially supported by NSF grant
DMS-9870027}

\email{nigel@math.missouri.edu}

\author{G. Lancien}
\address{Equipe de Math\'ematiques - UMR 6623, Universit\'e de
Franche-Comt\'e, F-25030 Besan\c con cedex}

\email{glancien@math.univ-fcomte.fr}


\begin{abstract}
We investigate the problem of $L^p$-maximal regularity on Banach
spaces having a Schauder basis. Our results improve those of a
recent paper.
\end{abstract}

\maketitle

\section{Introduction}

We will only recall the basic facts and definitions on maximal
regularity. For further information, we refer the reader to
\cite{COU}, \cite{DO}, \cite{LE1} or \cite{KL}.

We consider the following Cauchy problem:
$$\left\{
\begin{array}{ll}
u'(t)+B(u(t))=f(t)\ \ \ \ {\rm for\ } 0\leq t<T \\ u(0)=0
\end{array}\right.$$
where $T \in (0,+\infty)$, $-B$ is the infinitesimal generator of
a bounded analytic semigroup on a complex Banach space $X$ and $u$
and $f$ are $X$-valued functions on $[0,T)$. Suppose $1<p<\infty.$
$B$ is said to satisfy {\it $L^p-$maximal regularity} if whenever
$f\in L^p([0,T);X)$ then the solution $$ u(t) =\int_0^t
e^{-(t-s)B}f(s)\,ds$$ satisfies  $u'\in L^p([0,T);X).$  It is
known that $B$ has $L^p$-maximal regularity for some $1<p<\infty$
if and only if it has $L^p$-maximal regularity for every
$1<p<\infty$ \cite{DE}, \cite{DO}, \cite{SO}.  We thus say simply
that $B$ satisfies {\it maximal regularity (MR).}

As in \cite{KL}, we define:

\begin{Def} A complex Banach space $X$ has the {\it maximal
regularity property} (MRP) if $B$ satisfies (MR) whenever $-B$ is
the generator of a bounded analytic semigroup.
\end{Def}

Let us recall that De Simon \cite{DE} proved that any Hilbert
space has (MRP), and that the question whether $L^q$ for $1<q\neq
2<\infty$ has (MRP) remained open until recently. Indeed, in
\cite{KL} it is shown that a Banach space with an unconditional
basis (or more generally a separable Banach lattice) has (MRP) if
and only if it is isomorphic to a Hilbert space.

In this paper we attempt to work without these unconditionality
assumptions and study the (MRP) on Banach spaces with a
finite-dimensional Schauder
decomposition. In particular, we show that a UMD Banach space with an
(FDD)
 and satisfying (MRP) must be isomorphic to an
$\ell_2$ sum of finite dimensional spaces.

In the last question we consider the question of whether the solution $u$
of our Cauchy problem satisfies $u'\in L^2([0,T;L^r)$ if $f\in L^2([0,T);L^s)$.

This work was done during a visit of the  second author to the
Department of Mathematics of the University of Missouri in
Columbia in fall 1999; he would like to thank the Department for
its warm hospitality.

\section{Notation and background}

We will follow the notation of \cite{KL}. Let us now introduce
more precisely a few notions.

 If $F$ is a subset of the Banach space $X$, we denote by
$[F]$ the closed linear span of $F$.
 We denote by $(\eps_k)_{k=0}^\infty$ the standard
sequence of Rademacher functions on $[0,1]$ and by $(h_k)_{k=0}^\infty$
the standard Haar functions on $[0,1]$ (for convenience we index
from $0$).

 Let $1\leq p<\infty$. A Banach space $X$ has {\it type}
$p$ if there is a constant $C>0$ such that for every finite
sequence $(x_k)_{k=1}^K$ in $X$: $$(\int_0^1 \Vert
\sum_{k=1}^K\eps_k(t)x_k\Vert^2\, dt)^{1/2}\leq
C(\sum_{k=1}^K\Vert x_k\Vert^p)^{1/p}.$$ Notice that every Banach
space is of type 1.   A Banach space $X$ is called (UMD) if martingale
difference sequences in $L_2([0,1];X)$ are unconditional i.e. there is a
constant $K$ so that for every martingale difference sequence
$(f_n)_{n=1}^N$ we have
$$ \|\sum_{k=1}^N\delta_k f_k\|_{L_2(X)} \le K
\|\sum_{k=1}^Nf_k\|_{L_2(X)}$$ if $\sup_{k\le N}|\delta_k| \le 1.$

Let $(E_n)_{n\geq1}$ be a sequence of closed subspaces
of $X$. Assume that $(E_n)_{n\geq1}$ is a Schauder decomposition
of $X$ and let $(P_n)_{n\geq1}$ be the associated sequence of
projections from $X$ onto $E_n$. For convenience  we will also
denote this Schauder decomposition by $(E_n,P_n)_{n\geq1}$. The
decomposition constant is defined by $\sup_n\|\sum_{K=1}^nP_k\|$; this
is necessarily finite.  If each $(E_n)$ is finite-dimensional we refer to
$(E_n)$ as an (FDD) (finite-dimensional decomposition); an unconditional
(FDD) is abbreviated to (UFDD).

 If
$(E_n)_{n\ge 1}$ is a Schauder decomposition of $X$ and $(u_n)_{n=1}^N$
is a finite or infinite
sequence (i.e. $N\le \infty$) of the form
$u_n=\sum_{k=r_{n-1}+1}^{r_n}x_k$ where
$x_k\in E_k$ and $1=r_0<r_1<..<r_n<..$,
 then $(u_n)_{n\ge 1}$ is called a {\it block basic
sequence} of the decomposition $(E_n).$

 We denote by $\tree$ the set of all finite sequences of
positive integers, including the empty sequence denoted
$\emptyset$. For $a=(a_1,..,a_n)\in \tree$, $\vert a\vert=n$ is
the {\it length} of $a$ ($\vert \emptyset\vert=0$). For
$a=(a_1,..,a_k)$ (respectively $a=\emptyset$), we denote
$(a,n)=(a_1,..,a_k,n)$  (respectively $(a,n)=(n)$). A subset
$\beta$ of $\tree$ is a {\it branch} of $\tree$ if there exists
$(\sigma_n)_{n=1}^\infty \subset \Ndb$ such that
$\beta=\{(\sigma_1,..,\sigma_n);\ n\geq 1\}$.
 In this paper, for a Banach space $X$, we call a {\it
tree} in $X$ any family $(y_a)_{a\in \tree}\subset X$. A tree
$(y_a)_{a \in \tree}$ is {\it weakly null} if for any $a\in
\tree$, $(y_{(a,n)})_{n\geq 1}$ is a weakly null sequence.

 Let $(y_a)_{a \in \tree}$ be a tree in the Banach space
$X$. Let $T\subset \tree$, $(y_a)_{a \in T}$ is a {\it full
subtree} of $(y_a)_{a \in \tree}$ if $\emptyset \in T$ and for all
$a\in T$, there are infinitely many $n\in \Ndb$ such that
$(a,n)\in T$. Notice that if $(y_a)_{a \in T}$ is a full subtree
of a weakly null  tree $(y_a)_{a \in \tree}$, then it can be
reindexed as a weakly null  tree $(z_a)_{a \in \tree}$

\medskip
We now state a result of \cite{KL} that will be an
essential tool for this paper:

\begin{Thm}\label{sch} Let $(E_n,P_n)_{n \geq 1}$ be a Schauder
decomposition of the Banach space $X$. Let $Z_n=P_n^*X^*$ and
$Z=[\cup_{n=1}^{\infty}Z_n].$   Assume $X$ has (MRP). Then there
is a constant $C>0$ so that whenever $(u_n)_{n=1}^N$ are such that
$u_n\in [E_{2n-1},E_{2n}]$ and $(u_n^*)_{n=1}^N$ are such that
$u_n^*\in [Z_{2n-1},Z_{2n}]$ then $$
\left(\int_0^{2\pi}\|\sum_{n=1}^NP_{2n}u_n
e^{i2^nt}\|^2\frac{dt}{2\pi}\right)^{1/2}\le C\left( \int_0^{2\pi}
\| \sum_{n=1}^Nu_ne^{i2^nt}\|^2\frac{dt}{2\pi}\right)^{1/2}$$ and
$$ \left(\int_0^{2\pi}\|\sum_{n=1}^NP^*_{2n}u^*_n
e^{i2^nt}\|^2\frac{dt}{2\pi}\right)^{1/2}\le C\left( \int_0^{2\pi}
\| \sum_{n=1}^Nu^*_ne^{i2^nt}\|^2\frac{dt}{2\pi}\right)^{1/2}.$$
\end{Thm}
 We observe that, by a well-known result of Pisier
\cite{Pi1} these inequalities can be replaced by equivalent inequalities
(with a modified constant)
using $\eps_k$ in place of
$e^{i2^kt}:$
\begin{equation}\label{ineq1}
\|\sum_{n=1}^N P_{2n}u_n\eps_n\|_{L_2(X)}\le C \|\sum_{n=1}^N
u_n\eps_n\|_{L_2(X)}\end{equation} and
\begin{equation}\label{ineq2}
\|\sum_{n=1}^N P_{2n}u_n^*\eps_n\|_{L_2(X)}\le C \|\sum_{n=1}^N
u_n^*\eps_n\|_{L_2(X^*)}.\end{equation}

We refer the reader to \cite{W} for further recent developments in this
area.

\section{The main results}

We begin with a general result on spaces with a Schauder
decomposition:

\begin{Thm}\label{type} Let $X$ be a Banach space of type $p>1$ and with a
Schauder decomposition $(E_n)_{n=1}^\infty$. If $X$ has (MRP), then
there is a constant $C>0$ so that for any  block basic
sequence $(u_k)_{k=1}^N$ with respect to the decomposition $(E_n)$:
\begin{equation}\label{basic}\frac{1}{C}\sum_{k=1}^N \Vert u_k\Vert^2
\leq
\int_0^1
\Vert
\sum_{k=1}^N \eps_k(t)u_k\Vert^2\, dt \leq C\sum_{k=1}^K \Vert
u_k\Vert^2.\end{equation}
\end{Thm}

\begin{proof}
  If the result is false we can clearly inductively
construct an infinite normalized block basic sequence
$(u_n)_{n=1}^{\infty}$ so that there is no constant $C$ so that for
all finitely nonzero sequences $(a_k)_{k=1}^{\infty}$ we have:
\begin{equation}\label{ell2}
\frac{1}{C}\sum_{k=1}^N |a_k|^2 \leq \int_0^1 \Vert
\sum_{k=1}^N a_k\eps_k(t)u_k\Vert^2\, dt \leq C\sum_{k=1}^N
|a_k|^2\end{equation}
It therefore suffices to show that (\ref{ell2}) holds for every
normalized block basic sequence $(u_n)_{n=1}^{\infty}.$   We can clearly
then suppose $u_n\in E_n$.

We next use a theorem of Figiel and Tomczak-Jaegermann \cite{FT}
combined with
\cite{Pi2}
(see also
\cite{MS} p.112) that, since
$X$ has nontrivial
type for every $n\in\mathbb N$ there exists $\varphi(n)\in\mathbb N$ so
that any subspace $F$ of $X$ with dimension $\varphi(n)$ has a subspace
$H$ of dimension $n$ which is 2-complemented in $X$ and $2$-isomorphic to
$\ell_2^n.$

Assume (\ref{ell2}) is false. Then we can inductively find
a sequence $(a_n)_{n\ge 1}$ and an increasing sequence
$(r_n)_{n\ge 0}$ with $r_0=0$ so that
$r_{2n}>r_{2n-1}+\varphi(r_{2n-1}-r_{2n-2})$ for
$n\ge
1$,
$$\sum_{r_{2n}+1}^{r_{2n+1}}|a_k|^2 =1$$ and either
$$\int_0^1 \Vert
\sum_{k=r_{2n}+1}^{r_{2n+1}} a_k\eps_k(t)u_k\Vert^2\, dt >2^n $$ or
$$\int_0^1 \Vert
\sum_{k=r_{2n}+1}^{r_{2n+1}} a_k\eps_k(t)u_k\Vert^2\, dt <2^{-n}. $$

In order to create new Schauder decompositions of $X$, we will
need the following elementary lemma, that we state without a
proof:

\begin{Lem}\label{proj} Let $(E_n)_{n\geq 1}$ be a Schauder decomposition of a
Banach space $X$. Assume that each $E_n$ has a finite Schauder
decomposition $(F_{n,k})_{k=1}^{m_n}$ with a uniform bound on the
decomposition constant.  Then $(F_{1,1},\ldots,
F_{1,m_1},F_{2,1},\ldots,F_{2,m_2},\ldots)$ is also a Schauder
decomposition of $X.$
\end{Lem}

We denote the induced decomposition by
$\sum_{n=1}^{\infty}\oplus (\sum_{k=1}^{m_n}\oplus F_{n,k}).$

Now by assumption $E_{r_{2n-1}+1}+\cdots+E_{r_{2n}}$ which has dimension
at least $\varphi(r_{2n}-r_{2n-1})$ contains a subspace $H_n$ which is
$2$-Hilbertian and $2$-complemented in $X.$ Let $G_n$ be the complement
of $H_n$ in $E_{r_{2n-1}+1}+\cdots+E_{r_{2n}}$ by the projection of norm
2.  At the same time
$[u_k]$ is
1-complemented (by the Hahn-Banach theorem) in $E_k$ for $r_{2n-1}+1\le
k\le r_{2n}$  and let $F_k$ be its associated complement.  We thus have a
new Schauder decomposition:
$$(F_1,[u_1],F_2,[u_2],\ldots,F_{r_1},[u_{r_1}],H_1,G_1,F_{r_2+1},
[u_{r_2+1}],\ldots,[u_{r_3}],H_2,G_2,\ldots). $$   If we write
$D_n=F_{r_{2n-2}+1}+\cdots+F_{r_{2n-1}}+G_n$ then we have a Schauder
decomposition $$ \sum_{n=1}^{\infty}\oplus(D_n\oplus H_n\oplus
\sum_{k=r_{2n-2}+1}^{r_{2n}}\oplus [u_k]).$$

Next select a normalized basis $(v_k)_{k=r_{2n-2}+1}^{r_{2n-1}}$ of $H_k$
which is $2$-equivalent to the canonical basis of
$\ell_2^{r_{2n}-r_{2n-1}}.$  It is easy to see that we can obtain a  new
Schauder decomposition by interlacing the $(v_k)$ with the $(u_k)$ i.e.:
\begin{equation}\label{dec0} \sum_{n=1}^{\infty} (D_n\oplus
[u_{r_{2n-2}+1}]\oplus
[v_{r_{2n-2}+1}]\oplus\cdots\oplus [u_{r_{2n-1}}]\oplus
[v_{r_{2n-1}}]).\end{equation}
Now again using Lemma \ref{proj} we can form two further decompositions:
\begin{equation}\label{dec1} \sum_{n=1}^{\infty}( D_n\oplus
[u_{r_{2n-2}+1}+v_{r_{2n-2}+1}]\oplus
[v_{r_{2n-2}+1}]\oplus\cdots\oplus [u_{r_{2n-1}}+v_{r_{2n-1}}]\oplus
[v_{r_{2n-1}}]),\end{equation}
and
\begin{equation}\label{dec2} \sum_{n=1}^{\infty} (D_n\oplus
[u_{r_{2n-2}+1}+v_{r_{2n-2}+1}]\oplus
[u_{r_{2n-2}+1}]\oplus\cdots\oplus [u_{r_{2n-1}}+v_{r_{2n-1}}]\oplus
[u_{r_{2n-1}}]).\end{equation}

Now we can apply Theorem \ref{sch}.  If we use decomposition (\ref{dec1})
we note that $u_k=(u_k+v_k)-v_k$ and so for a suitable $C$ and all $n,$
$$\|\sum_{k=r_{2n-2}+1}^{r_{2n}}a_k(u_k+v_k)
\eps_k\|_{L_2(X)}\le C
\|
\sum_{k={r_{2n-2}+1}}^{r_{2n-1}}a_ku_k\eps_k\|_{L_2(X)}.$$
However, using decomposition (\ref{dec0}) there is also a constant $C'$
so that
$$\|\sum_{k=r_{2n-2}+1}^{r_{2n}}a_kv_k
\eps_k\|_{L_2(X)} \le C'
\|\sum_{k=r_{2n-2}+1}^{r_{2n}}a_k(u_k+v_k)
\eps_k\|_{L_2(X)}.$$  This leads to an
estimate:
$$ \left(\sum_{k={r_{2n-2}+1}}^{r_{2n-1}}|a_k|^2\right)^{\frac12}\le C_1
\|
\sum_{k={r_{2n-2}+1}}^{r_{2n-1}}a_ku_k\eps_k\|_{L_2(X)}.$$

If we use decomposition (\ref{dec2}) instead we obtain an estimate:
$$
\|
\sum_{k={r_{2n-2}+1}}^{r_{2n-1}}a_ku_ke^{i2^kt}\|_{L_2(X)}\le C_2
\left(\sum_{k={r_{2n-2}+1}}^{r_{2n-1}}|a_k|^2\right)^{\frac12}.$$

Combining gives us (\ref{ell2}) and completes the proof.\end{proof}

Let us first use this result to give a mild improvement of a result from
\cite {KL}:

\begin{Thm}\label{improv} Let $X$ be a reflexive space with an (FDD) and
with non-trivial type which embeds
into a space $Y$ with a (UFDD).  If $X$ has (MRP) then $X$ is isomorphic
to an $\ell_2-$sum of finite-dimensional spaces
$(\sum_{n=1}^{\infty}\oplus E_n)_{\ell_2}.$\end{Thm}

\begin{proof}
Using Proposition
1.g.4 of
\cite{LT}
(cf.
\cite{JZ})
we can block the given (FDD) to produce an (FDD) $(E_n)$ so that
$(E_{2n})_{n=1}^{\infty}$ and
$(E_{2n-1})_{n=1}^{\infty}$ are both (UFDD)'s.
Let us denote, as in Theorem \ref{sch}, the dual
(FDD) of $X^*$ by $(Z_n)_{n=1}^{\infty}.$
Now it follows applying
Theorem \ref{type} to both $X$ and $X^*$ (which also has (MRP)) that
there exists a constant
$C$ so that if
$x_n\in E_n$ and $x_n^*\in Z_n$ are two finitely nonzero sequences
\begin{align*}
 \|\sum_{k=1}^{\infty}x_{2k-j}\| &\le
C(\sum_{k=1}^{\infty}\|x_{2k-j}\|^2)^{\frac12}\\
 \|\sum_{k=1}^{\infty}x^*_{2k-j}\| &\le
C(\sum_{k=1}^{\infty}\|x^*_{2k-j}\|^2)^{\frac12}\end{align*} for $j=0,1.$
Hence
\begin{align*}
 \|\sum_{k=1}^{\infty}x_{k}\| &\le
2C(\sum_{k=1}^{\infty}\|x_{k}\|^2)^{\frac12}\\
 \|\sum_{k=1}^{\infty}x^*_{k}\| &\le
2C(\sum_{k=1}^{\infty}\|x^*_{k}\|^2)^{\frac12}.\end{align*}

Now for given $x_k$ we may find $y^*_k\in X^*$ with $\|y_k^*\|=\|x_k\|$
and $y_k(x_k^*)=\|x_k^*\|.$  Let $x_k^*=P_k^*y_k^*$ (where $P_k:X\to E_k$
is the projection associated with the FDD $(E_n)$).  Then $\|x_k^*\|\le C_1
\|x_k\|$ where $C_1 =\sup_n\|P_n\|<\infty.$  Hence if
$(x_k)_{k=1}^{\infty}$ is finitely nonzero, we have
$$ \|\sum_{k=1}^{\infty} x_k^*\| \le 2CC_1
(\sum_{k=1}^{\infty}\|x_k\|^2)^{\frac12}.$$  Thus
\begin{align*}
 \sum_{k=1}^{\infty}\|x_k\|^2 &= \sum_{k=1}^{\infty}x_k^*(x_k)\\
&= (\sum_{k=1}^{\infty}x_k^*)(\sum_{k=1}^{\infty}x_k)\\
&\le 2CC_1 (\sum_{k=1}^{\infty}\|x_k^*\|^2)^{\frac12}
\|\sum_{k=1}^{\infty}x_k\|\end{align*} so that we obtain the lower
estimate:
$$ (\sum_{k=1}^{\infty}\|x_k\|^2)^{\frac12} \le 2CC_1
\|\sum_{k=1}^{\infty} x_k\|.$$  This completes the proof.\end{proof}

We next give another application to (UMD)-spaces with (MRP).

\begin{Thm}\label{umd} Let $X$ be a (UMD) Banach space with an (FDD)
 satisfying (MRP). Then $X$ is isomorphic to an $\ell_2$-sum of
finite dimensional spaces, $(\sum_{n=1}^{\infty}\oplus E_n)_{\ell_2}.$
\end{Thm}

\begin{proof} Let $(E_n)$ be the given (FDD) of $X.$  We will show first
that there is a blocking $(F_n)$ of $(E_n)$ which satisfies an upper
$2-$estimate i.e. if there is a constant $A$ so that if $(x_n)$ is
block basic with respect to  $(F_n)$ and finitely non-zero then
\begin{equation}\label{upper}
\|\sum_{n=1}^{\infty}x_n\|\le A
(\sum_{n=1}^{\infty}\|x_n\|^2)^{\frac12}. \end{equation}  Once this is
done, the proof can be completed easily.  Indeed if $(Z_n)$ is the dual
decomposition to $(F_n)$ for $X^*$ then we can apply the fact that $X^*$
also has (MRP) ($X$ is reflexive) to block $(Z_n)$ to obtain a
decomposition which also has an upper $2$-estimate.  Thus we can assume
$(F_n)$ and $(Z_n)$ both have an upper 2-estimate and then repeat the
argument used in Theorem \ref{improv} to deduce that
$X=(\sum_{n=1}^{\infty}\oplus F_n)_{\ell_2}.$

Since $X$ necessarily has type $p>1$, we can apply Theorem \ref{type} and
assume $(E_n)$ obeys (\ref{basic}).

We now introduce a particular type of tree in the space $L_2([0,1);X).$
Let $\mathcal D_n$ for $n\ge 0$ be the sub-algebra of the Borel sets of
$[0,1)$ generated by the dyadic intervals $[(k-1)2^{-n},k2^{-n})$ for
$1\le k\le 2^n.$  Let $\mathbb E_n$ denote the conditional expectation
operator $\mathbb E_nf= \mathbb E(f|\mathcal D_n).$

We will say that a tree $(f_a)_{a\in\tree}$ is a {\it martingale
difference tree} or (MDT) if \begin{itemize} \item each $f_a$ is
$\mathcal D_{|a|}-$ measurable, \item
if $|a|>0$ then $\mathbb E_{|a|-1}f_a=0,$\item there exists $N$ so
that if
$|a|>N$ then $f_a=0.$\end{itemize}  In such a tree the partial sums along
any branch form a dyadic martingale which is eventually constant.

We will prove the following Lemma:

\begin{Lem}\label{tree} There is a constant $K$ so that if
$(f_a)_{a\in\tree}$ is a weakly null
(MDT),
there is a full subtree $(f_a)_{a\in T}$ so that for
any branch $\beta$ we have:
$$  \|\sum_{a\in\beta}f_a\|_{L_2(X)} \le
K(\sum_{a\in\beta}\|f_a\|_{L_2(X)}^2)^{\frac12}.$$ \end{Lem}

\begin{proof} For each $a$ we define integers $m_-(a)$ and $m_+(a)$.  If
$f_a\neq 0$ we set $m_-(a)$ to be the greatest $m$ so that
$$ \|\sum_{k=1}^mP_mf_a\|_{L_2(X)}\le 2^{-|a|-1}\|f_a\|_{L_2(X)}$$ and
$m_+(a)$ to be the least $m> m_-(a)$ so that
 $$ \|\sum_{k=m+1}^{\infty}P_kf_a\|_{L_2(X)}\le
2^{-|a|-1}\|f_a\|_{L_2(X)}.$$ If
$f_{\emptyset}=0$ we set $m_-(\emptyset)=0$ and $m_+(\emptyset)=1$; if
$f_a=0$
where $a\neq \emptyset$ we set $m_-(a)$ to be the last member of $a$ and
$m_+(a)=m_-(a)+1.$

Since $(f_a)$ is weakly null we have $\lim_{n\to\infty}m_-(a,n)=\infty$
for every $a.$  It is then easy to pick a full subtree $T$ so that
$m_-(a,n)>m_+(a)$ whenever $a, (a,n)\in T.$  Now let
$g_a=\sum_{k=m_-(a)+1}^{m_+(a)}f_a.$  Then $\|f_a-g_a\|_{L_2(X)} \le
2^{-|a|}\|f_a\|_{L_2(X)}.$

For any branch $\beta$ of $T$, we have that $g_a(t)$ is a block basic
sequence with respect to $(E_n)$ for every $0\le t<1.$  Hence
$$ \left(\int_0^1 \|\sum_{a\in\beta}\epsilon_{|a|}(s)g_a(t)\|_X^2
ds\right)^{\frac12} \le C
\left(\sum_{a\in\beta}\|g_a(t)\|_X^2\right)^{\frac12}.$$  Integrating
again we have
$$  \left (\int_0^1
\|\sum_{a\in\beta}\epsilon_{|a|}(s)g_a\|_{L_2(X)}^2ds\right)^{\frac12}
\le C \left(\sum_{a\in\beta}\|g_a\|_{L_2(X)}^2\right)^{\frac12}.$$
From this we get
$$  \left (\int_0^1 \|
\sum_{a\in\beta}\epsilon_{|a|}(s)f_a\|_{L_2(X)}^2ds\right)^{\frac12}
\le 2C
\left(\sum_{a\in\beta}\|f_a\|_{L_2(X)}^2\right)^{\frac12}+\sum_{a\in\beta}
2^{-|a|}\|f_a\|.$$

Estimating the last term by the Cauchy-Schwarz inequality and using the
fact that $X$ is (UMD) we get the Lemma.\end{proof}

Now we introduce a functional $\Phi$ on $X$ by defining
$\Phi(x)$ to be the infimum of all $\lambda>0$ so that for every weakly
null
(MDT) $(f_a)_{a\in\tree}$ with $f_{\emptyset}=x\chi_{[0,1)}$ we have  a full subtree
$T$ so that for any branch $\beta$
\begin{equation}\label{Fdef}
  \|\sum_{a\in\beta}f_a\|^2_{L_2(X)} \le  \lambda  + 2K^2
\sum_{\substack{a\in\beta\\ a\neq
\emptyset}}\|f_a\|^2_{L_2(X)}.\end{equation}

Note that since
 $$ \|\sum_{a\in\beta}f_a\|^2_{L_2(X)} \le 2(\|x\|^2+
 \|\sum_{\substack{a\in\beta\\ a\neq\emptyset}}f_a\|^2_{L_2(X)})$$ we
have
an estimate $\Phi(x)\le 2\|x\|^2$.  By considering the null tree we have
$F(x)\ge \|x\|^2.$  It is clear that $\Phi$ is continuous and
$2$-homogeneous.  Most importantly we observe that $\Phi$ is convex; the
proof of this is quite elementary and we omit it.
It follows that we can define an equivalent norm by $\tnorm{x}^2=\Phi(x)$
and $\|x\|\le \tnorm{x}\le 2\|x\|$ for $x\in X.$

Next we prove that if $x\in X$ and $(y_n)$ is a weakly null sequence then
\begin{equation}\label{weaknull}
\limsup_{n\to\infty}(\tnorm{x+y_n}^2+\tnorm{x-y_n}^2) \le
2\tnorm{x}^2+4K^2\limsup_{n\to\infty}\|y_n\|^2.\end{equation}

We first note that we can suppose $\lim_{n\to\infty}\tnorm{x\pm y_n}$
and
$\lim_{n\to\infty}\|y_n\|^2$ all exist. Now suppose $\epsilon>0.$
Then we can find weakly null
(MDT)'s
$(f^n_a)_{a\in \tree}$ with $f^n_{\emptyset}\equiv x+y_n$ so that for every full
subtree $T$ we have a branch $\beta$ on which:
\begin{equation}\label{negFdef}
\|\sum_{a\in\beta}f^n_a\|_{L_2(X)}^2+\epsilon
>
 \tnorm{x+y_n}^2+2K^2
\sum_{\substack{a\in\beta\\
a\neq\emptyset}}\|f^n_a\|_{L_2(X)}^2.\end{equation}
In fact by easy induction we can pick a full subtree so that
(\ref{negFdef}) holds for {\it every} branch.  Hence we suppose the
original tree satisfies (\ref{negFdef}) for every branch.

Similarly we may find weakly null (MDT)'s $(g^n_a)_{a\in\tree}$ with
$g^n_{\emptyset}\equiv x-y_n$ and for every branch $\beta,$
\begin{equation*} \|\sum_{a\in\beta}g^n_a\|_{L_2(X)}^2
+\epsilon> \tnorm{x-y_n}^2+2K^2
\sum_{\substack{a\in\beta\\
a\neq\emptyset}}\|g^n_a\|_{L_2(X)}^2.\end{equation*}

We next consider the (MDT) defined by
$h_{\emptyset}\equiv x$,
$$ h_{(n)}(t) =\begin{cases} y_n\qquad \text{if } 0\le t<\frac12\\ -y_n
\qquad\text{if
}
\frac12\le t<1\end{cases}$$  and if $|a|>1$ then
$$ h_{(a,n)}(t) =\begin{cases} f^n_a(2t-1) \qquad\text{if } 0\le
t<\frac12\\ g^n_a(2t) \qquad\text{if } \frac12\le t<1. \end{cases} $$

Now for every branch of the (MDT) $(h_a)_{a\in\tree}$ with initial
element $\{n\}$ we have
$$ \|\sum_{a\in\beta}h_a\|_{L_2(X)}^2 + \epsilon >
\frac12(\tnorm{x+y_n}^2+\tnorm{x-y_n}^2) + 2K^2
\sum_{\substack{a\in\beta\\
|a|>1}}\|h_a\|_{L_2(X)}^2.$$

However, from the definition of $\Phi(x)=\tnorm{x}^2$ it follows that
there exists
$n_0$
so that if $n\ge n_0$ we can find a branch $\beta$ whose initial element
is $n$ and such that
$$ \|\sum_{a\in\beta}h_a\|_{L_2(X)}^2< \tnorm{x}^2+  2K^2
\sum_{\substack{a\in\beta\\
|a|>0}}\|h_a\|_{L_2(X)}^2 +\epsilon.$$   Combining gives the equation
(for $n\ge n_0)$,
$$ \frac12(\tnorm{x+y_n}^2+\tnorm{x-y_n}^2) \le \tnorm{x}^2
+2K^2\|y_n\|^2+2\epsilon.$$ This proves (\ref{weaknull}).
But note that if
$y_n$ is weakly null we have $\liminf_{n\to\infty} \tnorm{x-y_n} \ge
\tnorm{x}$ and so we deduce:
$$\limsup_{n\to\infty}\tnorm{x+y_n}^2 \le \tnorm{x}^2
+4K^2\limsup\|y_n\|^2.$$

Using this equation it is now easy to block the Schauder decomposition
$(E_n)$  to produce a Schauder decomposition $(F_n)$ with the property
that for any $N$ if $x\in F_1+\cdots+F_N$ and $y\in
\sum_{k=N+2}^{\infty}F_k$ then
$$ \tnorm{x+y} \le (1+\delta_N)(\tnorm{x}^2+4K^2\|y\|^2)^{\frac12},$$
where
$\delta_N>0$ are chosen to be decreasing and  so that
$\prod_{N=1}^{\infty}(1+\delta_N)\le
2.$
Next suppose $(x_k)$ is any finitely non-zero block basic sequence with
respect to $(F_n).$   By an easy induction we obtain for
$j=0,1$:
$$ \tnorm{\sum_{k=1}^nx_{2k-j}} \le 4K^2
\prod_{k=1}^{n-1}(1+\delta_{2k-j})(\sum_{k=1}^n\tnorm{x_{2k-j}}^2)^{\frac12}.$$
Hence
$$  \|\sum_{k=1}^nx_k\| \le 32K^2 (\sum_{k=1}^n\|x_k\|^2)^{\frac12}.$$
 This establishes (\ref{upper}) and as shown earlier this suffices to
complete the proof.\end{proof}

{\it Remark.}  Recently Odell and Schlumprecht \cite{OS} showed that a
separable
Banach space $X$ can be embedded in an $\ell_p-$sum of finite-dimensional
spaces for $1<p<\infty$ if and only if $X$ is reflexive and every
normalized weakly null tree has a branch which is equivalent to the usual
$\ell_p-$basis.  This result is closely related to the proof of
the previous theorem.

\section{On $L^r$-regularity in $L^s$ spaces}

Let $s\in [1,\infty)$. We consider our usual Cauchy problem:

$$\left\{
\begin{array}{ll}
u'(t)+B(u(t))=f(t)\ \ \ \ {\rm for\ } 0\leq t<T \\ u(0)=0
\end{array}\right.$$
where $T \in (0,+\infty)$, $-B$ is the infinitesimal generator of
a bounded analytic semigroup on $L^s=L^s([0,1])$ and $f\in L^2([0,T);L^s)$.
Then we ask the following question: for what values of $s$ and $r$ in $[1,\infty)$
does the solution $$ u(t) =\int_0^t
e^{-(t-s)B}f(s)\,ds$$ necessarily satisfies  $u'\in L^p([0,T);L^r)$? Thus we
introduce the following definition:

\begin{Def} Let $r$ and $s$ in $[1,\infty)$. We say that $(r,s)$ is a {\it regularity pair}
if whenever $-B$ is the infinitesimal generator of
a bounded analytic semigroup on $L^s=L^s([0,1])$ and $f\in L^2([0,T);L^s)$, the solution $u$
of $$\left\{
\begin{array}{ll}
u'(t)+B(u(t))=f(t)\ \ \ \ {\rm for\ } 0\leq t<T \\ u(0)=0
\end{array}\right.$$
satisfies $u'\in L^p([0,T);L^r)$.
\end{Def}

Notice that it follows from previous results (\cite{DE}, \cite{LE1} and \cite{KL})
that $(s,s)$ is a regularity pair if and only if $s=2$. This is extended by our next result:

\begin{Thm}\label{reg}   Let $r$ and $s$ in $[1,\infty)$. Then $(r,s)$ is
a regularity pair if and only if $r\leq  s=2$.
\end{Thm}

\begin{proof} It follows clearly from the work of De Simon \cite{DE}, that if $r\leq  s=2$
then $(r,s)$ is a regularity pair.

\smallskip\noindent So let now $(r,s)$ be a regularity pair. Since $L^1$ does not have (MRP)
(\cite{LE1}), we have that $s>1$.  Then, solving our Cauchy problem with $B=0$, we obtain that
$r\leq s$. Thus we can limit ourselves to the case $s>1$ and $1\leq r\leq s$.

\noindent Then by the closed graph Theorem, for any $B$ so that
$-B$ is the infinitesimal generator of a bounded analytic
semigroup on $L^s=L^s([0,1])$, there is a constant $C>0$ such that
for any $f\in L^2([0,T);L^s)$: $$\Vert u'\Vert_{L^2(L^s)} \leq
C\Vert f\Vert_{L^2(L^s)}.$$

Using the inclusion $L^s\subset L^r$ for $r\leq s$, we can now state the following
analogue of Theorem 2.1:

\begin{Prop}\label{pair} Let $(E_n,P_n)_{n \geq 1}$ be a Schauder
decomposition of $L^s$. Assume that $(r,s)$ is a regularity pair. Then there
is a constant $C>0$ so that whenever $(u_n)_{n=1}^N$ are such that
$u_n\in [E_{2n-1},E_{2n}]$  then $$
\|\sum_{n=1}^N P_{2n}u_n\eps_n\|_{L^2(L^r)}\le C \|\sum_{n=1}^N
u_n\eps_n\|_{L^2(L^s)}$$
\end{Prop}

Then our first step will be to show that the Haar system satisfies some lower-2 estimates
in $L^s$ in the following sense:

\begin{Lem}\label{low} If there exists $r\leq s$ such that $(r,s)$ is  a
regularity pair, and if $s<p<2$ or $p=2$ then there is
a constant $C>0$ such that for any normalized block basic sequence $(v_1,...,v_n)$ of $(h_k)$
and for any $a_1,..,a_n$ in $\Cdb$:
$$\Vert \sum_{k=1}^n a_kv_k\Vert_{L^s} \geq C(\sum_{k=1}^n \vert
a_k\vert^p)^{\frac1p}.$$
\end{Lem}

\begin{proof}  We first observe that if $1<p<2,$
it follows from the work of J. Bretagnolle, D.
Dacunha-Castelle and J.L. Krivine \cite{BCDK} on $p$-stable random
variables that there is a  sequence $(e_n)_{n\geq 1}$ in $L^1$
which is equivalent to the canonical basis of $\ell_p$ in any
$L^q$ for $1\leq q<p$. Thus $(e_n)$ is weakly null in $L^s$, and
by a gliding hump argument, we may assume that $(e_n)$ is actually a
block basic sequence with respect to the Haar basis.  If $p=2$ then the
Rademacher functions already form a block basic sequence in every $L^q$
for $1\le q<\infty.$

 Now assume the Lemma is false. We pick a normalized
block basic sequence
$(v_1,...,v_{n_1})$ of $(h_k)$
and $a_1,..,a_{n_1}$ in $\Cdb$ so that
$$\Vert \sum_{k=1}^{n_1} a_kv_k\Vert_{L^s} \leq (\sum_{k=1}^{n_1} \vert
a_k\vert^p)^{\frac1p}=1.$$
Then pick $m_1\in \Ndb$ such that $(v_1,..,v_{n_1},e_{m_1})$ is a
block basic sequence of $(h_k)$.
 By induction, we pick a normalized block basic sequence
$(v_{n_j+1},...,v_{n_{j+1}})$ of $(h_k)$,
$a_{n_j+1},..,a_{n_{j+1}}$ in $\Cdb$ and $m_{j+1}\in \Ndb$ so that
$(v_1,..,v_{n_1},\eps_{m_1},v_{n_1+1},..,v_{n_{j+1}},\eps_{m_{j+1}})$
is a block basic sequence of $(h_k)$ and $$\Vert
\sum_{k=n_j+1}^{n_{j+1}} a_kv_k\Vert_{L^s} \leq
\frac{1}{2^j}(\sum_{k=n_j+1}^{n_{j+1}} \vert
a_k\vert^p)^{\frac1p}=\frac{1}{2^j}.$$ So we can find $(I_k)_{k\geq 1}$
and
$(J_k)_{k\geq 1}$ two sequences of finite intervals of $\Ndb$ such
that $\{I_k,J_k:\ k\geq 1\}$ is a partition of $\Ndb$ and for all
$k\geq 1$, $v_k\in [h_j,\ j\in I_k]$ and $e_{m_k}\in [h_j,\
j\in J_k]$. Then set $$X_k=[h_j:~ j\in I_k\cup J_k].$$ Then
$(X_k)$ is an unconditional Schauder decomposition of $L^s$. Each
$X_k$ can be decomposed into $X_k=E_{2k-1}\oplus E_{2k}$, where
$E_{2k-1}=[v_k+\eps_{m_k}]$, $e_{m_k}\in E_{2k}$ and the
corresponding projections are uniformly bounded. So, by Lemma
\ref{proj}, $(E_k)_{k\geq 1}$ is a Schauder decomposition of
$L^s$. We can now make use of Proposition \ref{pair}. If we
decompose $a_kv_k=a_k(v_k+e_{m_k})-a_ke_{m_k}$ in
$E_{2k-1}\oplus E_{2k}$, we obtain that there is a constant $C>0$
such that for all $n\geq 1$: $$\Vert \sum_{k=1}^n
a_kv_k\eps_k\Vert_{L^2(L^s)} \geq C(\sum_{k=1}^n \vert
a_k\vert^p)^{\frac1p}.$$ Since $(v_k)$ is an unconditional basic
sequence in $L^s$, there is a constant $K>0$ so that for all
$n\geq 1$: $$\Vert \sum_{k=1}^n a_kv_k\Vert_{L^s} \geq
K(\sum_{k=1}^n \vert a_k\vert^p)^{\frac1p},$$ which is in contradiction
with our construction.
\end{proof}

We now conclude the proof of Theorem \ref{reg}.  The Haar basis of $L^s$
has a block basic sequence equivalent to the standard basis of
$\ell_{\max(s,2)}.$  Hence Lemma \ref{low} shows that $\max(s,2)\le p$
whenever $s<p<2$ or $p=2.$  Thus $s=2.$
\end{proof}

\end{document}